\documentclass[11pt,a4paper]{article}

\usepackage{psfrag}
\usepackage{amsmath,amssymb,cite}
\usepackage{latexsym}
\usepackage{graphicx}
\usepackage{lscape}
\usepackage{afterpage}

\usepackage{hyperref}
\hypersetup{
    colorlinks=true,
    linkcolor=blue,
    filecolor=magenta,      
    urlcolor=cyan,
    citecolor=red,
}
\DeclareMathOperator*{\argmin}{argmin}
\DeclareMathOperator*{\tv}{TV}
\DeclareMathOperator*{\sgn}{sgn}

\newcommand{\ds}{\displaystyle}
\newcommand{\nexto}{\kern -0.54em}
\newcommand{\dR}{{\rm {I\ \nexto R}}}

\newcommand{\dZ}{{\cal Z \kern -0.7em Z}}
\newcommand{\dC}{{\rm\hbox{C \kern-0.8em\raise0.2ex\hbox{\vrule
height5.4pt width0.7pt}}}}
\newcommand{\dQ}{{\rm\hbox{Q \kern-0.85em\raise0.25ex\hbox{\vrule
height5.4pt width0.7pt}}}}
\newcommand{\proofbox}{\hspace{\fill}{$\Box$}}
\newtheorem{lemma}{Lemma}
\newtheorem{theorem}{Theorem}

\newtheorem{remark}{Remark}

\newenvironment{proof}{Proof.}{\proofbox}

\begin{document}

\author{Authors}

\author{
C. Yal{\c c}{\i}n Kaya\footnote{School of Information Technology and Mathematical Sciences, University of South Australia, Mawson Lakes, S.A. 5095, Australia. E-mail: yalcin.kaya@unisa.edu.au\,.}
}

\title{\vspace{-10mm}\bf Optimal Control of the Double Integrator with Minimum Total Variation}

\maketitle

\begin{abstract} {\noindent\sf We study the well-known minimum-energy control of the double integrator, along with the simultaneous minimization of the total variation in the control variable. We derive the optimality conditions and obtain the unique optimal solution to the combined problem, where the initial and terminal boundary points are specified.  We study the problem from a multi-objective optimal control viewpoint, constructing the Pareto front.   We show that the unique asymptotic optimal control function, for the minimization of the total variation alone, is piecewise constant with one switching at the midpoint of the time horizon.  For any instance of the boundary conditions of the problem, we prove that the asymptotic optimal total variation is exactly $2/3$ of the total variation of the minimum-energy control.  We illustrate the results for a particular instance of the problem and include a link to a video which animates the solutions while moving along the Pareto front.}
\end{abstract}

\begin{verse}
\end{verse}

\begin{verse} 
\end{verse}

\pagestyle{myheadings}
\markboth{}{\sf\scriptsize Optimal Control of the Double Integrator with Minimum Total Variation\ \ by C. Y. Kaya}

\section{Introduction}

\begin{figure}[t]
\begin{center}
\psfrag{u}{$u(t)$}
\psfrag{x10}{$x_1(0) = s_0$}
\psfrag{x1}{$x_1(t) := y(t)$}
\psfrag{x2}{$x_2(t) := \dot{y}(t)$}
\includegraphics[width=100mm]{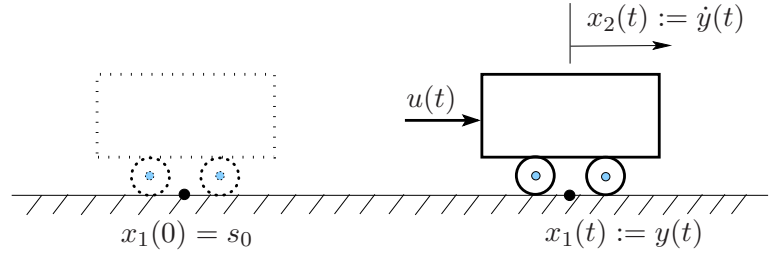} 
\end{center}
\caption{A  simplified physical model of a car as a point mass.}
\label{fig:car}
\end{figure}

The double integrator is a mathematical model for a point mass, typically idealizing a car in rectilinear motion on a flat and frictionless plane as schematically illustrated in Figure~\ref{fig:car}.  It also constitutes a model for analogous rotational-mechanical and electrical systems~\cite{Wellstead2000}.  One should recall that a cubic curve between two oriented points, which minimizes its averaged acceleration, or more precisely, the $L^2$-norm of its acceleration, serves as a building block for cubic splines~\cite{Dontchev1993, OpfObe1988}.  This latter case can be represented as the energy-minimizing double integrator.

Due to its simplicity, optimal control of the double integrator is studied virtually in every course of lectures on optimal control theory.  In the teaching of optimal control theory and its applications, although the minimum-energy, minimum-effort and minimum-time control of the double integrator are widely studied, minimization of total variation is not even considered, presumably because a maximum principle for the control minimizing its total variation does not exist.

The double integrator model is so simple that an analytical solution can be worked out easily for the problem of energy minimization.  Moreover, for the case of minimum-time control, where the control variable is bound-constrained, the optimal control can simply be shown to be bang--bang with at most one switching, i.e., the control variable switches from one bound to the other, and it does so at most once.  The control structure can also be worked out easily in the case of minimum-effort control, where the $L^1$-norm of the control function is minimized.  In summary, optimal control of the double integrator yields simple but rich-enough examples for illustrations of some key aspects of the theory of optimal control~\cite{RaoBer2001}.

Total variation of a function can be broadly described as the total vertical distance traversed by the graph of the function (a precise definition is to be given in Section~\ref{total_variation}).  A small total variation in the control function is obviously desirable, as it would make the control system easier to design and implement, resulting in, for example, smaller or lighter motors for a robot or a spacecraft.

Although there is a lack of theory and results for the pure minimization of total variation, it is often imposed in addition to the minimization of another functional, for instance, energy or duration of time.  This is done in the earlier works \cite{TeoJen1991, LoxLinTeo2013,WanYuTeo2016}, where the optimal control problem is discretized directly by assuming piecewise-constant optimal control variables.  This discretization simplifies the expression for the total variation in control; however, optimality conditions for the original (continuous-time) problem cannot be derived or verified, because of the discretization itself.

Total variation is widely used as a regularization term in more general optimization problems such as imaging and signal processing (see \cite{GinTanShiHenZha2019} and the references therein).  It has also relatively recently been used as a regularization term for parameter estimation in linear quadratic control \cite{KosKur2013}. A bound on the total variation in the control is derived for minimum-time linear control problems in \cite{Silin1982}, although the total variation itself is not incorporated into the minimization problem.

In the present article, in addition to the minimization of energy, we consider the minimization of the total variation in the control variable of the double integrator.  In other words, we aim to study simultaneous minimization of energy and total variation, giving rise to {\em multi-objective optimization} and the study of the set of all trade-off/compromise solutions called the {\em Pareto front}.  Optimal control problems which involve total variation have not been studied yet from the viewpoint of multi-objective optimal control.

In this paper, we use a tutorial approach.  First, in Section~\ref{sec:minen_control}, we introduce the double integrator model as well as the problem of energy minimization as an optimal control problem.  This is a standard problem in optimal control; so, we derive the optimal solution without going into details.  

In Section~\ref{problem}, we define the total variation of a function and state the energy and total variation minimization problem, by appending the total variation in control as a weighted term to the energy functional.  Next, we augment the state variable vector, so that the problem can be rewritten and posed as an optimal control problem in standard form.  We derive optimality conditions,  and discuss the problem as a multi-objective optimal control problem.    By means of asymptotic analysis, we derive an optimal solution for the pure total variation minimization problem. These kinds of results on total variation do not exist in the literature.

In Section~\ref{sec:example}, via an example instance of the problem, we illustrate the results given in Section~\ref{problem}.
In particular, we provide, via a URL link in \cite{video}, a video illustration of the multi-objective solutions on the Pareto front, so that evolution of the solutions as the weight of total variation is varied can be animated and observed.

Finally, in Section~\ref{sec:conclusion}, we offer concluding remarks and provide various relevant open problems.

\section{Minimum-Energy Control}
\label{sec:minen_control}

Consider the car as a point unit mass, moving on a frictionless planar ground in a fixed line of action, as shown in Figure~\ref{fig:car}.  Let the position of the car at time $t$ be given by $y(t)$ and the velocity by $\dot{y}(t):=(dy/dt)(t)$.  By Newton's second law of motion, $\ddot{y}(t) = u(t)$, where $u(t)$ is the summation of all the external forces applied on the car, in this case the force simply representing the acceleration and deceleration of the car.   This differential equation model is referred to as the {\em double integrator} in system theory literature, since $y(t)$ can be obtained by integrating $u(t)$ twice.

Let $x_1 := y$ and $x_2 := \dot{y}$.  The problem of minimizing the energy of the car, which starts at a position $x_1(0) = s_0$ with a velocity $x_2(0) = v_0$ and finishes at the final position $x_1(1) = s_f$ with velocity $x_2(1) = v_f$, within one unit of time, can be posed as follows.
\[
\mbox{(Pe) }\left\{\begin{array}{rl}
\ds\min & \ \ \ds\frac{1}{2}\int_0^1 u^2(t)\,dt  \\[5mm] 
\mbox{subject to} & \ \ \dot{x}_1(t) = x_2(t)\,,\ \ x_1(0) = s_0\,,\ \ x_1(1) = s_f\,, \\[2mm]
& \ \ \dot{x}_2(t) = u(t)\,,\ \ \ \,x_2(0) = v_0\,,\ \ x_2(1) = v_f\,.
\end{array} \right.
\]
Here, the functions $x_1$ and $x_2$ are referred to as the {\em state variables} and $u$ the {\em control variable}. As a first step in writing the conditions of optimality for this optimization problem, define the Hamiltonian function $H$ for Problem (Pe) in the usual way as
\begin{equation}  \label{Hamiltonian0}
H(x_1,x_2,u,\lambda_1,\lambda_2) := \frac{1}{2}\,u^2 + \lambda_1\,x_2 + \lambda_2\,u\,,
\end{equation}
where $\lambda(t) := (\lambda_1(t),\lambda_2(t))\in\dR^2$ is the {\em adjoint variable} (or {\em costate}) {\em vector} such that (see~\cite{Hestenes66})
\begin{equation} \label{adjoint} 
\dot{\lambda}_1 = -\partial H /\partial x_1\quad\mbox{and}\quad
\dot{\lambda}_2 = -\partial H /\partial x_2\,.
\end{equation} 
The equations in~\eqref{adjoint} simply reduce to 
\begin{equation} \label{adjoint_sol} 
\lambda_1(t) = \bar\lambda_1\quad\mbox{and}\quad \lambda_2(t) = -\bar\lambda_1\,t - c\,, 
\end{equation} 
where $\bar\lambda_1$ and $c$ are real constants.  By calculus of variations, or the maximum principle with an unconstrained control variable (see~\cite{Hestenes66}), if $u$ is optimal, then
\begin{equation}  \label{opt_en}
\partial H /\partial u = 0\,,\mbox{ i.e.,}\quad u(t) = -\lambda_2(t) = \bar\lambda_1\,t + c\,.
\end{equation}
Substituting $u(t)$ in \eqref{opt_en} into the differential equations and solving these equations by also utilizing the boundary conditions in Problem~(Pe), one gets the analytical solution
\begin{eqnarray}
u(t) &=& \bar\lambda_1\,t + c\,,  \label{en_u} \\[1mm]
x_1(t) &=& \frac{1}{6}\,\bar\lambda_1\,t^3 + \frac{1}{2}\,c\,t^2 + v_0\,t + s_0\,,  \label{en_x1} \\[1mm]
x_2(t) &=& \frac{1}{2}\,\bar\lambda_1\,t^2 + c\,t + v_0\,,  \label{en_x2}
\end{eqnarray} 
for all $t\in[0,1]$, where
\begin{eqnarray}
\bar\lambda_1 &=& -12\,(s_f - s_0) + 6\,(v_0 + v_f)\,,  \label{en_4} \\[1mm]
c &=& 6\,(s_f - s_0) - 2\,(2\,v_0 + v_f)\,.  \label{en_5}
\end{eqnarray}

We note that the position variable $x_1(t)$ of the car is a cubic polynomial of time.  Therefore, the minimum-energy control solution, despite being so simple, constitutes a building block for the problem of finding a cubic spline interpolant passing through a given set of points.

\section{Minimization of Total Variation}
\label{problem}

\subsection{Total variation of a function}
\label{total_variation}

\begin{figure}[t]
\begin{center}
\psfrag{u}{$u(t)$}
\psfrag{t}{$t$}
\includegraphics[width=100mm]{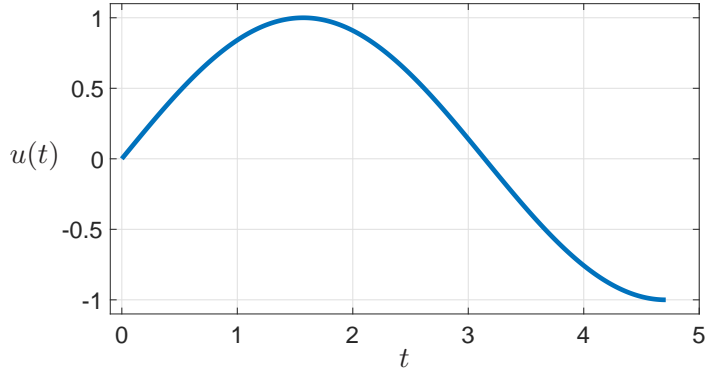} 
\end{center}
\caption{Graph of $u(t) = \sin t$ over $[0,3\pi/2]$, illustrating that $\tv(u) = 3$.}
\label{sint}
\end{figure}

The {\em total variation} of a function $u:[t_0,t_f]\to\dR$ is defined as
\begin{equation} \label{totalvar1}
\tv(u) := \sup\sum_{i=1}^N |u(t_i) - u(t_{i-1})|\,,
\end{equation}
where the supremum is taken over all partitions
\begin{equation} \label{partition}
t_0 < t_1 < \cdots < t_N = t_f
\end{equation}
of the interval $[t_0,t_f]$ (see \cite{Kreyszig1978}).  Here, $N\in\{1,2,3,\ldots\}$ is arbitrary as is the choice of the values $t_1,\cdots,t_{N-1}$ in $[t_0,t_f]$ which, however, must satisfy \eqref{partition}.  The function $u$ is said to be of {\em bounded variation} on $[t_0,t_f]$, if $\tv(u)$ is finite.  If $u$ is absolutely continuous on $[t_0,t_f]$, in other words, $u\in W^{1,1}([t_0,t_f];\dR)$, then
\begin{equation} \label{totalvar2}
\tv(u) = \int_{t_0}^{t_f} |\dot{u}(t)|\,dt\,,
\end{equation}
where $\dot{u} := du/dt$.  Practically speaking, $\tv(u)$ as given in \eqref{totalvar2} represents the total distance traversed by the projection of the $u(t)$ vs.\ $t$ graph along the vertical $u(t)$ axis.  Figure~\ref{sint} illustrates this interpretation with $u(t) = \sin t$ over $[0,3\pi/2]$, where clearly $\tv(u) = 3$.

\newpage

\subsection{Minimum-Total-Variation Control of the Double Integrator}

Recall that in the case when one has Problem~(Pe), minimizing only the energy, the solution is as given in~\eqref{en_u}--\eqref{en_5}.  So, clearly $\tv(u) = 6\,|2\,(s_f - s_0) - v_f - v_0|$.

We consider optimal control problems where we aim to minimize the total variation in the control variables {\em in addition to} the energy functional.
\[
\mbox{(Ptv) }\left\{\begin{array}{rl}
\ds\min & \ \ \ds\frac{1}{2}\int_0^1 u^2(t)\,dt  + \alpha\,\tv(u) \\[5mm] 
\mbox{subject to} & \ \ \dot{x}_1(t) = x_2(t)\,,\ \ x_1(0) = s_0\,,\ \ x_1(1) = s_f\,, \\[2mm]
& \ \ \dot{x}_2(t) = u(t)\,,\ \ \ \,x_2(0) = v_0\,,\ \ x_2(1) = v_f\,,
\end{array} \right.
\]
where $\alpha > 0$ is referred to as the {\em weight}.  We assume that $u$ is absolutely continuous on $[0,1]$, in other words, $u\in W^{1,1}([0,1])$.  Then we define the new control variable $v(t) := \dot{u}(t)$ for a.e. $t\in[0,1]$.  Using \eqref{totalvar2}, Problem~(Ptv) can now be reformulated by incorporating the new variable as
\[
\mbox{(Paug) }\left\{\begin{array}{rl}
\ds\min & \ \ \ds\frac{1}{2}\int_0^1 \left(u^2(t)\,dt  + \alpha\,|v(t)| \right) dt\\[5mm] 
\mbox{subject to} & \ \ \dot{x}_1(t) = x_2(t)\,,\ \ x_1(0) = s_0\,,\ \ x_1(1) = s_f\,, \\[2mm]
& \ \ \dot{x}_2(t) = u(t)\,,\ \ \ \,x_2(0) = v_0\,,\ \ x_2(1) = v_f\,,\\[2mm]
& \ \ \dot{u}(t) = v(t)\,,\mbox{\ \ for a.e. } t\in[0,t_f]\,.
\end{array} \right.
\]
In this augmented form of the problem, $u$ becomes a new state variable.

\subsection{Optimality Conditions}

The {\em Hamiltonian} function for Problem~(Paug) is given by
\begin{equation} \label{Hamiltonian}
H(x_1,x_2,u,v,\lambda_1,\lambda_2,\eta) := \frac{1}{2}\,u^2 + \alpha\,|v| + \lambda_1\,x_2 + \lambda_2\,u + \eta\,v\,,
\end{equation}
where $\lambda(t) = (\lambda_1(t),\lambda_2(t))\in\dR^2$ and $\eta(t)\in\dR$ are {\em adjoint variables} defined by (see~\cite{Hestenes66})
\begin{eqnarray}
&& \dot{\lambda}_1 := -\partial H /\partial x_1 = 0\quad\mbox{and}\quad \dot{\lambda}_2 := -\partial H /\partial x_2 = -\lambda_1\,,  \label{lambda1a} \\[1mm]
&& \ \,\dot{\eta} := -\partial H /\partial u \ \,= -u - \lambda_2\,,  \quad \eta(0) = 0\,,\ \ \eta(1) = 0\,,  \label{etaa} 
\end{eqnarray}
In other words,
\begin{eqnarray}
&& \lambda_1(t) = \bar{\lambda}_1\,,\quad\mbox{and}\quad \lambda_2(t) = -\bar{\lambda}_1\,t - c\,, \label{lambd1b} \\[1mm]
&& \ \,\dot{\eta}(t) = -u(t) + \bar{\lambda}_1\,t + c\,,  \quad \eta(0) = 0\,,\ \ \eta(1) = 0\,.  \label{etab} 
\end{eqnarray}
where $\bar\lambda_1$ and $c$ are real constants.  Note that, although the expressions in \eqref{lambda1a} are respectively the same as those in \eqref{adjoint}, the real constants $\bar\lambda_1$ and $c$ in this case depend on the value of $\alpha$ and so are different in general.

Next we state the maximum principle (see~\cite[Theorem 1.5.1]{Vinter2000}) for our setting as follows.  

\noindent {\bf Maximum Principle.}\ \ Suppose that $x_1,x_2,u\in W^{1,1}([0,1];\dR)$ and $v\in L^1([0,1];\dR)$ solve Problem~(Paug).  Then there exist functions $\lambda_1,\lambda_2,\eta\in W^{1,1}(0,t_f;\dR)$ such that\linebreak $(\lambda_1(t), \lambda_2(t), \eta(t)) \neq \bf0$, for every $t\in[0,1]$, and, in addition to the state differential equations and other constraints given in Problem~(Paug) and the adjoint differential equations in~\eqref{lambda1a}--\eqref{etaa}, the following condition holds:
\begin{equation}  \label{optcont1}
v = \argmin_{w\in\dR}\ H(x_1,x_2,u,w,\lambda_1,\lambda_2,\eta) = \argmin_{w\in\dR}\ (\alpha\,|w| + \eta\,w)\,;\mbox{\ \ for a.e. } t\in[0,1]\,.
\end{equation}
Condition~\eqref{optcont1} implies that
\begin{equation}  \label{optcont2}
v(t) = \left\{\begin{array}{rl}
0\,, &\ \ \mbox{if\ \ } |\eta(t)| < \alpha\,, \\[2mm]
\mbox{undetermined}\,, &\ \ \mbox{if\ \ } |\eta(t)| = \alpha\,,
\end{array} \right.
\end{equation}
for a.e. $t\in[0,1]$.  Note that $|\eta_i(t)| > \alpha$ is not allowed by the maximum principle, as otherwise one would get $v(t) = -\infty$.  

In view of \eqref{optcont2}, when $-\alpha < \eta(t) < \alpha$,\ \ a.e.\ \ $t\in[0,1]$, the original control $u(t)$ is (possibly piecewise) constant.  What if $|\eta(t)| \equiv \alpha$ over a subinterval of $[0,1]$?  If so, then we refer to the optimal control in this subinterval as singular control, which we elaborate further next.

\noindent {\bf Singular control.}\ \ If there exist $s_1$ and $s_2$ such that $|\eta(t)| = \alpha$ for every $t\in[s_1,s_2] \subset [0,t_f]$ (in fact, one has either $\eta(t) = \alpha$ or $\eta(t) = -\alpha$ for every $t\in[s_1,s_2]$, because of the continuity of $\eta$), then the control variable $v(t)$ for every $t\in[s_1,s_2]$ is said to be {\em singular}.  A candidate for a {\em singular optimal control} $v(t)$ might be obtained by observing that, since $\eta(t)$ is constant over $[s_1,s_2]$, one will have $\dot\eta(t) = \ddot\eta(t) = 0$ for every $t\in[s_1,s_2]$.  By using \eqref{etab}, this observation yields
\[
\dot\eta(t)\equiv 0 = -u(t) - \lambda_2(t)\,,
\]
i.e.,
\[
u(t) = \bar\lambda_1\,t + c\,,
\]
and so
\[
v(t) = \bar\lambda_1\,,
\]
for all $t\in[s_1,s_2]$.

\noindent {\bf Optimal control.}\ \ With the incorporation of the singular control, and by the continuity of the adjoint variable~$\eta$, \eqref{optcont2} can be rewritten as
\begin{equation}  \label{optcont_convex1a}
v(t) = \left\{\begin{array}{rl}
0\,, &\ \ \mbox{if\ \ } |\eta(t)| < \alpha\,, \\[2mm]
\bar\lambda_1\,, &\ \ \mbox{if\ \ } |\eta(t)| = \alpha\,,
\end{array} \right.
\end{equation}
for all $t\in[0,1]$.  Note that $v(t)$ in \eqref{optcont_convex1a} is piecewise-constant and so $u(t)$ is piecewise-linear and continuous in $t$.  Then, by \eqref{etab}, $\eta(t)$ is continuous and piecewise-quadratic in $t$.  Note in particular that, differentiating both sides of the ODE in \eqref{etab}, using $\dot{u} = v$ and substituting
\eqref{optcont_convex1a}, one gets
\begin{equation}  \label{eta_ddot}
\ddot{\eta}(t) = \left\{\begin{array}{rl}
\bar\lambda_1\,, &\ \ \mbox{if\ \ } |\eta(t)| < \alpha\,, \\[2mm]
0\,, &\ \ \mbox{if\ \ } |\eta(t)| = \alpha\,.
\end{array} \right.
\end{equation}
The expression in \eqref{eta_ddot} and the boundary conditions in \eqref{etab} imply that there will be at most two junction points, $0 < t_1 < t_2 < 1$, for $\eta(t)$.  Namely, either $\eta(t) = \alpha$ or $\eta(t) = -\alpha$, for $t_1\le t < t_2$, and $\eta(t)$ is quadratic in $t$, for $0\le t < t_1$ and $t_2\le t \le 1$, with the
same constant second derivative $\bar\lambda_1$.  In other words,
\begin{equation}  \label{v_P1a}
v(t) = \left\{\begin{array}{rl}
0\,, &\ \ \mbox{if\ \ } 0\le t < t_1\ \ \mbox{or}\ \ t_2\le t \le 1\,, \\[2mm]
\bar\lambda_1\,, &\ \ \mbox{if\ \ } t_1\le t < t_2\,.
\end{array} \right.
\end{equation}
Then from $\dot{u} = v$ and continuity of $u$, one gets
\begin{equation}  \label{optcont_convex1}
u(t) = \left\{\begin{array}{ll}
\bar{u}_1\,, &\ \ \mbox{if\ \ } 0\le t < t_1\,, \\[2mm]
\bar{u}_1 + \bar\lambda_1\,(t-t_1)\,, &\ \ \mbox{if\ \ } t_1\le t < t_2\,, \\[2mm]
\bar{u}_3\,, &\ \ \mbox{if\ \ } t_2\le t \le 1\,.
\end{array} \right.
\end{equation}
where $\bar{u}_1$ and $\bar{u}_3$ are unknown constants.  Subsequently, $c = -\bar\lambda_1\,t_1 + \bar{u}_1$,
\[
\lambda_2(t) = \bar\lambda_1\,(t_1-t) - \bar{u}_1\,,
\]
and
\begin{equation}  \label{eta_convex1}
\eta(t) = \left\{\begin{array}{ll}
\ds\frac{1}{2}\bar\lambda_1\,(t^2 - 2\,t_1\,t)\,, &\ \ \mbox{if\ \ } 0\le
      t < t_1\,, \\[2mm] 
\alpha\mbox{\ \ or\ \ }-\alpha\,, &\ \ \mbox{if\ \ } t_1\le t \le t_2\,,
      \\[2mm] 
\ds\frac{1}{2}\bar\lambda_1\left[t^2 + 2\,t_2\,(1-t) - 1\right], &\ \
      \mbox{if\ \ } t_2 < t \le 1\,.
\end{array} \right.
\end{equation}
Note that $\lim_{t\to t_1^-}\eta(t) = \lim_{t\to t_2^+}\eta(t)$ (both equal to $\alpha$ or $-\alpha$), which, after simple algebraic manipulations, yields
\begin{equation} \label{junctions}
t_1 = 1 - t_2\,.
\end{equation}
We also note that $\lim_{t\to t_1^-}\eta(t) = -\bar\lambda_1\,t_1^2 / 2 = \mp\alpha$, i.e.,
\begin{equation}  \label{lambda1}
\bar\lambda_1 = \pm \frac{2\,\alpha}{t_1^2}\,.
\end{equation}

\begin{lemma}  \label{lem:t1}
One has that $0 < t_1 < 1/2$\,.
\end{lemma}
\begin{proof}
The proof is furnished by the fact that $0 < t_1 < t_2 < 1$ and \eqref{junctions}.
\end{proof}

\subsection{Multi-Objective Optimal Control}

Problem~(Ptv), or equivalently Problem~(Paug), concerns a {\em simultaneous} minimization of two objectives, which can simply be written as 
\begin{equation}  \label{MOP}
\mbox{(Pmo) }\ \ \ \min_{u\in\cal{U}}\ \ \left[ \varphi_1(u)\,,\ \varphi_2(u)\right]\,,
\end{equation}
where
\begin{equation}  \label{individuals}
\varphi_1(u) := \frac{1}{2}\int_0^1 u^2(t)\,dt\quad\mbox{and}\quad\varphi_2(u) := TV(u)\,.
\end{equation}
Problem~(Pmo) is referred to as a {\em multi-objective}, or {\em vector}, {\em optimal control problem}, with $\cal{U}$ representing the feasible, or admissible, set of all control functions satisfying the differential equation constraints and the boundary conditions---see \cite{BonKay2010} and the references therein.  The set of all solutions of~\eqref{MOP} is usually infinite, consisting of all trade-off, or Pareto, solutions.  Broadly speaking, a {\em Pareto} solution is a solution where one cannot improve the value of one objective functional without making the other worse.  The set of all Pareto solutions in the $\varphi_1\varphi_2$-plane (or the value space) is referred to as the {\em Pareto front} of Problem~(Pmo).  An example of a Pareto front is given in Figure~\ref{fig:convex}(a) (see details in Section~\ref{sec:example}).

For solving \eqref{MOP}, a typical approach is to consider a scalarization of the vector objective and so reduce Problem~(Pmo) to a single-objective optimal control problem.  Note that $\varphi_1$ and $\varphi_2$ are convex and the constraint set represents linear differential equations and linear boundary conditions.  Therefore we can use the weighted-sum scalarization (see \cite{BonKay2010}):
\[
\mbox{(Ps1) }\ \ \ \min_{u\in\cal{U}}\ \ \alpha_1\,\varphi_1(u) + (1-\alpha_1)\,\varphi_2(u)\,,
\]
where $\alpha_1\in(0,1)$.  Since $\alpha_1 \neq 0$, we can define $\alpha := (1-\alpha_1)/\alpha_1$ and write
\[
\mbox{(Ps2) }\ \ \ \min_{u\in\cal{U}}\ \ \varphi_1(u) + \alpha\,\varphi_2(u)\,,
\]
with $\alpha\in(0,\infty)$.  We note that Problems~(Ps1) and (Ps2) are equivalent and that Problem~(Ps2) is in the same form as Problem~(Ptv).

In this case, the individual functionals in \eqref{individuals} can be calculated using~\eqref{optcont_convex1} and \eqref{junctions}, in terms of the unknown parameters $t_1$, $\bar{u}_1$ and $\bar{u}_3$, as follows.
\begin{eqnarray}
&& \varphi_1(u)  = \frac{1}{2}\left[(\bar{u}_1^2+\bar{u}_3^2)\,t_1 + \frac{1}{3\,\bar\lambda_1}\,\left(\bar{u}_3^3 - \bar{u}_1^3\right)\right],  \label{phi1} \\[3mm]
&&  \varphi_2(u) = |\bar{u}_1 - \bar{u}_3|\,.  \label{phi2}
\end{eqnarray}

\subsection{Solution}
\label{solution}

By using~\eqref{optcont_convex1} and the initial conditions in Problem~(Ptv), and by integrating directly, one can obtain the following expressions for the state variables $x_1(t)$ and $x_2(t)$.
\begin{equation}  \label{x1}
x_1(t) = \left\{\begin{array}{ll}
\ds\frac{1}{2} \bar{u}_1\,t^2 + v_0\,t + s_0\,, &\ \ \mbox{if\ \ } 0\le t < t_1\,, \\[3mm] 
\ds\frac{1}{6}\bar\lambda_1(t-t_1)^3 + \frac{1}{2} \bar{u}_1\,t^2 + v_0\,t + s_0\,, &\ \ \mbox{if\ \ } t_1\le t < t_2\,, \\[2mm] 
\ds\frac{1}{6}\bar\lambda_1(t_2-t_1)^3 + \frac{1}{2} \bar{u}_1\,t_2^2 + v_0\,t_2 + s_0 & \\[2mm]
\ \ \ +\ \ds\frac{1}{2} \bar{u}_3(t-t_2)^2 + \left[\frac{1}{2}\bar\lambda_1(t_2-t_1)^2 + \bar{u}_1\,t_2 + v_0\right](t-t_2)\,, &\ \
      \mbox{if\ \ } t_2\le t \le 1\,;
\end{array} \right.
\end{equation}

\begin{equation}  \label{x2}
x_2(t) = \left\{\begin{array}{ll}
\bar{u}_1\,t + v_0\,, &\ \ \mbox{if\ \ } 0\le t < t_1\,, \\[3mm] 
\ds\frac{1}{2}\bar\lambda_1(t-t_1)^2 + \bar{u}_1\,t + v_0\,, &\ \ \mbox{if\ \ } t_1\le t < t_2\,, \\[2mm] 
\ds\frac{1}{2}\bar\lambda_1(t_2-t_1)^2 + \bar{u}_1\,t_2 + v_0 + \bar{u}_3(t-t_2)\,, &\ \
      \mbox{if\ \ } t_2\le t \le 1\,.
\end{array} \right.
\end{equation}

Finally, writing out the terminal conditions $x_1(1) = 0$ and $x_2(1) = 0$ using \eqref{x1} and \eqref{x2}, respectively, using \eqref{junctions} and \eqref{lambda1}, and carrying out lengthy manipulations, we obtain the following.
\begin{eqnarray}
&&    4\,t_1^3 - 3\,\left(2 \pm \frac{v_f + v_0 + 2\,(s_0 - s_f)}{\alpha}\right) t_1^2 + 1 = 0 \label{cubic} \\[1mm]
&&    \bar\lambda_1 = \pm\frac{2\,\alpha}{t_1^2}\,, \label{lambda1}   \\[1mm]
&&    \bar{u}_1 = v_f - v_0 - \frac{\bar\lambda_1}{2}(1-2\,t_1)\,, \label{u1}  \\[1mm]
&&    \bar{u}_3 = 2\,(v_f - v_0) - \bar{u}_1\,.  \label{u3}
\end{eqnarray}
Once $t_1$ is determined as a solution of~\eqref{cubic}, the parameters $\bar{\lambda}_1$, $\bar{u}_1$ and $\bar{u}_3$ in \eqref{lambda1}--\eqref{u3}, respectively, can explicitly be found.  The following lemma guides us as to which of the signs~$\pm$ (in the coefficient of the $t_1^2$-term) in~\eqref{cubic} will yield a solution and that whether the solution will be unique or not.

\begin{lemma}[Existence and uniqueness of the solution of a cubic equation]  \label{lem:cubic}
Let $c$ be a real constant.  Then the equation
\begin{equation}  \label{cubic1}
4\,t^3 - 3\,\left(2 + c \right) t^2 + 1 = 0
\end{equation}
has a unique solution for $c>0$, and has no solution for $c<0$, over the interval~$(0,1/2)$.
\end{lemma}
\begin{proof}
Let $f_c(t) := 4\,t^3 - 3\,\left(2 + c \right) t^2 + 1$.  Suppose that $c = 0$.  Then $f_0(0) = 1 > 0$, $f_0(1) = -1 < 0$ and that $f_0'(t) = 12\,t\,(\,t - 1) < 0$ for all $t\in(0,1)$, implying that $f(t) = 0$ for a unique $t\in(0,1)$.  One can easily check that $f_0(p_i) = 0$, $i = 1,2,3$, where $p_1 = (1-\sqrt{3})/2 < 0$, $p_2 = 1/2$, $p_3 = (1+\sqrt{3})/2 > 1$.  We have that $f_c'(t) = 12\,t\,(t - 1 - c/2)$ and that
\[
f_c'(t)\left\{\begin{array}{rl}
< f_0'(t) & \mbox{ if}\ \ c > 0\,, \\[2mm]
> f_0'(t) & \mbox{ if}\ \ c < 0\,,
\end{array}\right.
\]
for all $t\in(0,1)$. \\
(i) Suppose $c > 0$.  Then $f_c(0) = 1 > 0$ and $f_c(1/2) = 3/2 - 3\,(2+c)/4 = -3\,c/4 < 0$.  With $f_c'(t) < f_0'(t) < 0$ for all $t\in(0,1/2)$, we conclude that $f_c(t)$ has a unique zero in $(0,1/2)$. \\
(ii) Suppose $c < 0$.  Then $f_c(0) = 1 > 0$.  Since the only zero $f_0$ has is $1/2$ in the interval $[0,1/2]$, and $f_c'(t) > f_0'(t)$ for all $t\in(0,1/2)$,  $f_c$ has no zero in $[0,1/2]$.  
\end{proof}

\begin{remark}  \rm
Lemma~\ref{lem:t1} states that $t_1\in(0,1/2)$, and Lemma~\ref{lem:cubic} implies which sign in~\eqref{cubic1} needs to be considered in order to find a unique $t_1$.  By comparing \eqref{cubic1} and \eqref{cubic}, it is immediate to see that the sign of $(v_f + v_0 + 2\,(s_0 - s_f))$ has to be taken into account.  The following theorem provides the unique solution to Problem~(Ptv), based on this observation. 
\proofbox
\end{remark}

\begin{theorem}[Solution of Problem~(Ptv)]  \label{thm:solution}
The solution to Problem~(Paug) is unique and given by the expressions for the optimal control variable in~\eqref{optcont_convex1}, and the state variables in \eqref{x1}--\eqref{x2}, where the parameter $t_1$ is the solution of the cubic in~\eqref{cubic1} on the interval $(0,1/2)$, with $c = \bar{c}$ such that
\begin{equation}  \label{c}
\bar{c} = \frac{|v_f + v_0 + 2\,(s_0 - s_f)|}{\alpha}\,,
\end{equation}
the parameter $\bar\lambda_1$ given by
\begin{equation}  \label{lambda1_thm}
\bar\lambda_1 = \sgn(\bar{c})\,\frac{2\,\alpha}{t_1^2}\,,
\end{equation}
and the parameters $\bar{u}_1$ and $\bar{u}_3$ given by \eqref{u1}--\eqref{u3}.  As a result, the optimal total variation is given by
\[
\tv(u) = \left|v_0 - v_f  - \sgn(\bar{c})\,\frac{4\,\alpha}{t_1^2}\,(1 - 2\,t_1)\right|\,.
\]
\end{theorem}
\begin{proof}
By Lemma~\ref{lem:cubic} with $c = \bar{c} > 0$, there exist a unique $t_1\in(0,1/2)$ which solves \eqref{cubic1} and satisfies the optimality condition in Lemma~\ref{lem:t1}.  Recall again Lemma~\ref{lem:cubic} that, for $c < 0$, \eqref{cubic1} has no solution in $(0,1/2)$.  Therefore, in~\eqref{cubic}, we use the plus sign when $(v_f + v_0 + 2\,(s_0 - s_f)) > 0$, and the minus sign when $(v_f + v_0 + 2\,(s_0 - s_f)) < 0$.  Subsequently, this argument transforms \eqref{cubic} into \eqref{cubic1} with $c = \bar{c}$.  Furthermore, the $\pm$ sign in~\eqref{lambda1} is replaced by $\sgn(\bar{c})$ accordingly, yielding \eqref{lambda1_thm}.  The rest of the theorem follows from direct substitutions.
\end{proof}

\subsection{Asymptotic Solution (as \boldmath$\alpha \to\infty$)}
\label{sec:asymptotic}

As mentioned in the Introduction, it is not possible to write down the necessary conditions of optimality for the minimization of the total variation in the control variable alone.  Nevertheless, an analytic solution of Problem~(Ptv) can still be obtained by studying the asymptotic behaviour of the solutions when $\alpha\to\infty$.  In this case, Equation~\eqref{cubic} becomes $4\,t_1^3 - 6\,t_1^2 + 1 = 0$, which has three real roots: 1/2 and $(1\pm\sqrt{3})/2$.  This means that, in $(0,1/2)$, $t_1\to 1/2$.  Then, by \eqref{junctions}, $t_2\to1/2$.  Moreover, from Equation~\eqref{lambda1}, $\bar\lambda_1\to\pm\infty$.  However, these limit values of $t_1$ and $\bar\lambda_1$ make the expression in \eqref{u1} indeterminate.  Therefore, we need to write the asymptotic expressions for the state variables (with $t_1 = t_2 = 1/2$), in order to proceed:
\begin{equation}  \label{x1_asymp}
x_1(t) = \left\{\begin{array}{ll}
\ds\frac{1}{2} \bar{u}_1\,t^2 + v_0\,t + s_0\,, &\ \ \mbox{if\ \ } 0\le
      t < 1/2\,, \\[3mm] 
\ds\frac{1}{8}\,\bar{u}_1 + \frac{1}{2}\,v_0 + s_0 + \left(\frac{1}{2} \bar{u}_1 + v_0\right)\left(t - \frac{1}{2}\right)+ \frac{1}{2}\,
                  \bar{u}_3\left(t-\frac{1}{2}\right)^2\,, &\ \
      \mbox{if\ \ } 1/2\le t \le 1\,,
\end{array} \right.
\end{equation}

\begin{equation}  \label{x2_asymp}
x_2(t) = \left\{\begin{array}{ll}
\bar{u}_1\,t + v_0\,, &\ \ \mbox{if\ \ } 0\le t < 1/2\,, \\[3mm] 
\ds\frac{1}{2}\,\bar{u}_1 + v_0 + \bar{u}_3\left(t-\frac{1}{2}\right)\,, &\ \
      \mbox{if\ \ } 1/2\le t \le 1\,.
\end{array} \right.
\end{equation}

Now we can state the result, as $\alpha\to\infty$, in the following theorem.
\begin{theorem}[Asymptotic minimum total variation]   \label{thm:asymp_solution}
The unique asymptotic optimal control variable $u(t)$ of Problem~(Ptv), as $\alpha\to\infty$, is piecewise constant with a single switching at $t = 1/2$, namely
\begin{equation}  \label{u_asymp}
u(t) = \left\{\begin{array}{rl}
4\,(s_f - s_0) - v_f - 3\,v_0\,, &\ \ \mbox{if\ \ } 0\le t < 1/2\,, \\[1mm] 
3\,v_f + v_0 - 4\,(s_f - s_0)\,, &\ \ \mbox{if\ \ } 1/2\le t \le 1\,.
\end{array} \right.
\end{equation}
Consequently, the asymptotic optimal total variation is given by
\begin{equation}  \label{tv_asymp}
\tv(u) = 4 \left|2\,(s_f - s_0) - v_f - v_0\right|\,.
\end{equation}
\end{theorem}
\begin{proof}
The boundary conditions $x_1(1) = s_f$ and $x_2(1) = v_f$ using \eqref{x1_asymp}--\eqref{x2_asymp} yield, after manipulations,
\begin{eqnarray*}
3\,\bar{u}_1 + \bar{u}_3 &=& 8\,(s_f - s_0 - v_0)\,,  \\
\bar{u}_1 + \bar{u}_3 &=& 2\,(s_f - s_0 - v_0)\,,
\end{eqnarray*}
the solution of which is $\bar{u}_1 = 4\,(s_f - s_0) - v_f - 3\,v_0$ and $\bar{u}_3 = 3\,v_f + v_0 - 4\,(s_f - s_0)$, as required by \eqref{u_asymp}.  The switching time, $t = t_1 = 1/2$, is found as explained in the first paragraph of this subsection~\ref{sec:asymptotic}.  The expression in~\eqref{tv_asymp} is obtained by simply substituting the solutions for $\bar{u}_1$ and $\bar{u}_3$ above into $\tv = |\bar{u}_1 - \bar{u}_3|$.
\end{proof}

\begin{remark}  \label{rem:asymptotic} \rm
Recall that when no minimization of the total variation in control is done, i.e., when only the energy is minimized, the total variation is $\tv(u) = 6\,|2\,(s_f - s_0) - v_f - v_0|$.  It is interesting to note that the asymptotic minimum total variation in~\eqref{tv_asymp} is exactly $2/3$ of the total variation in minimum-energy control.
\proofbox
\end{remark}

\section{An Example}
\label{sec:example}

To demonstrate the results in Theorems~\ref{thm:solution} and \ref{thm:asymp_solution}, as well as illustrate what the Pareto front looks like using the expressions in \eqref{phi1}--\eqref{phi2}, we consider a particular instance when $s_0 = 0$, $s_f = 0$, $v_0 = 1$ and $v_f = 0$.  In view of the interpretation of the double integrator dynamics provided in the Introduction, this particular instance means that the car with an initial unit velocity is required to come to rest in the same position where it started the motion. 

The minimum energy solution can be obtained directly, after substituting $s_0 = 0$, $s_f = 0$, $v_0 = 1$ and $v_f = 0$ into \eqref{en_u}--\eqref{en_5}, as
\begin{eqnarray*}
u(t) &=& 6\,t - 4\,, \\[1mm]
x_1(t) &=& t^3 - 2\,t^2 + t\,, \\[1mm]
x_2(t) &=& 3\,t^2 - 4\,t + 1\,,
\end{eqnarray*}
for $t\in[0,1]$.  In this case, clearly, $\tv(u) = 6$\,.

\begin{figure}[t]
\begin{minipage}{80mm}
\begin{center}
\psfrag{a1}{\small $\alpha = 10^{-6}$}
\psfrag{a2}{\small $\alpha = 0.05$}
\psfrag{a3}{\small $\alpha = 0.4$}
\psfrag{a4}{\hspace*{1mm}\small $\alpha = 10^6$}
\psfrag{f1}{$\varphi_1$}
\psfrag{f2}{$\varphi_2$}
\includegraphics[width=80mm]{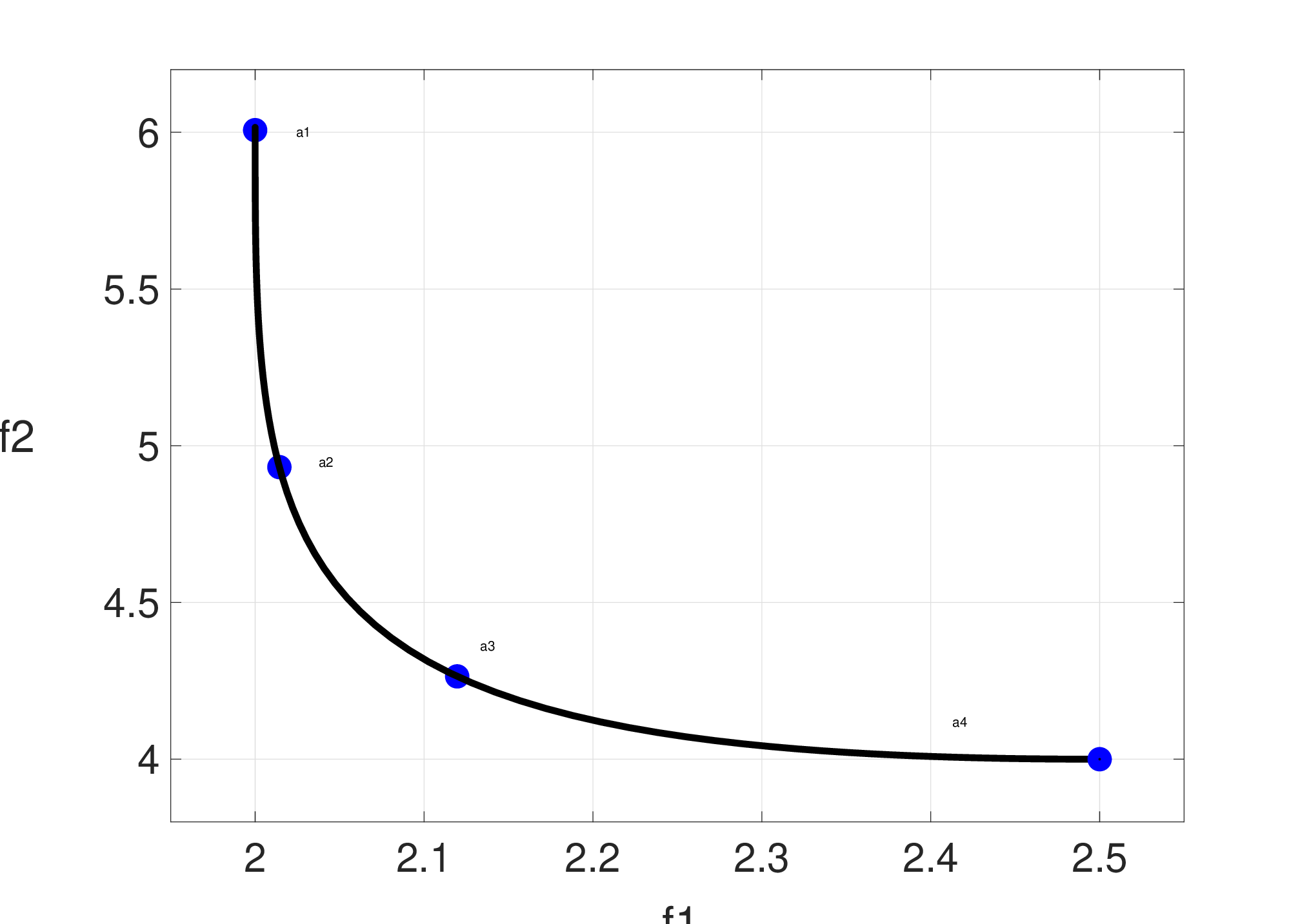} \\[3mm]
(a) The Pareto front
\end{center}
\end{minipage}
\hspace{-3mm}
\begin{minipage}{80mm}
\begin{center}
\psfrag{u}{\hspace*{-0.5mm}$u(t)$}
\psfrag{t}{$t$}
\psfrag{a0}{\small $\alpha = 10^{-6}$}
\psfrag{a5}{\small $\alpha = 10^{6}$}
\includegraphics[width=85mm]{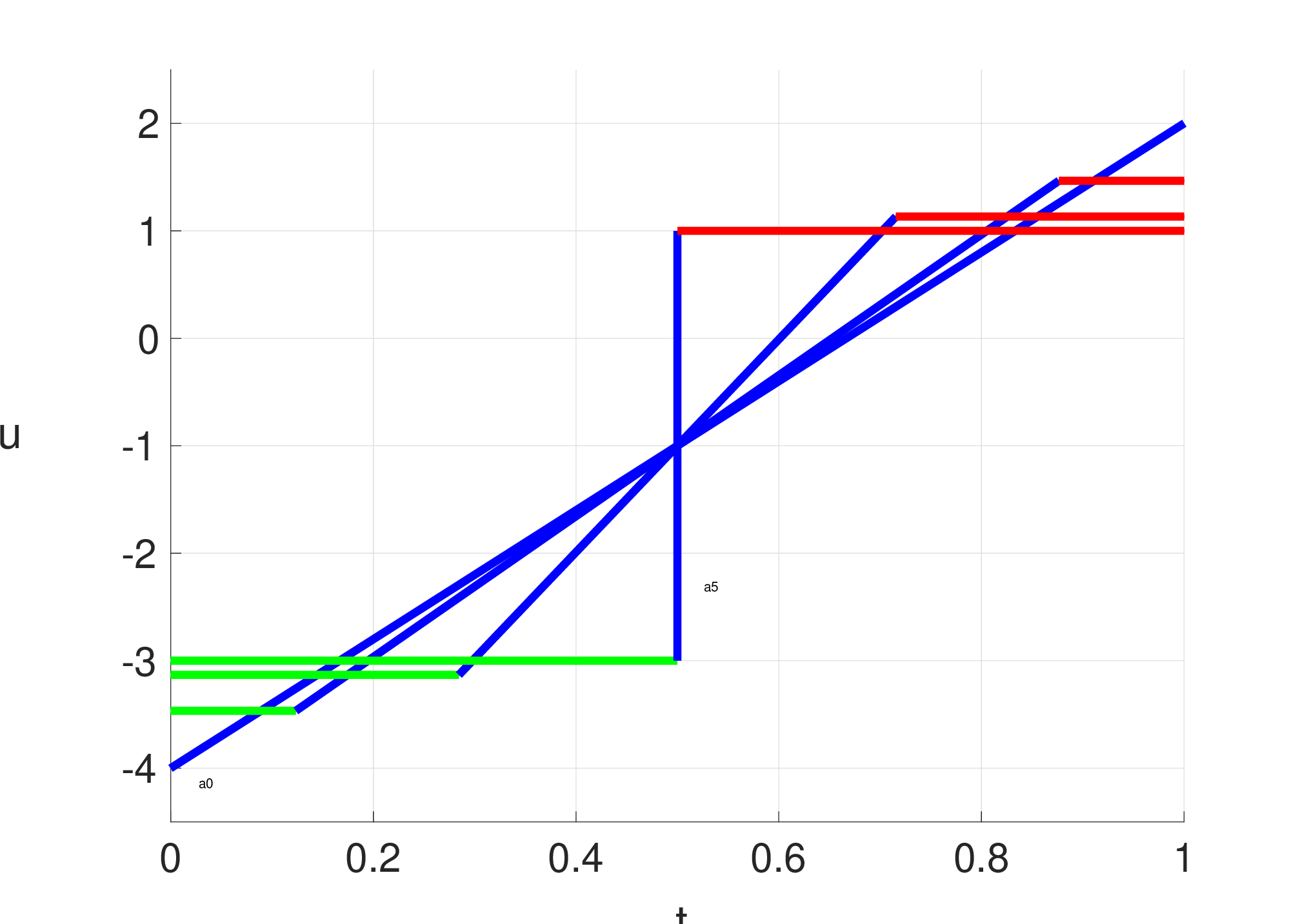} \\[3mm]
(b) The control variable
\end{center}
\end{minipage}
\
\caption{\sf The Pareto front and the control variable for the multi-objective problem, with $s_0 = 0$, $s_f = 0$, $v_0 = 1$ and $v_f = 0$.} 
\label{fig:convex}
\end{figure}

Figure~\ref{fig:convex} depicts the full Pareto front, as well as the optimal control variable for the parameter values $\alpha = 10^{-6}$, 0.05, 0.4 and $10^6$.  In drawing the graphs, first, Theorem~\ref{thm:solution} has been used: the unknown parameters $t_1$, $\bar\lambda_1$, $\bar{u}_1$, $\bar{u}_3$ (and $t_2 = 1 - t_1$).  Then $u(t)$, $\varphi_1(u)$ and $\varphi_2(u)$ have been computed as given in \eqref{optcont_convex1} and \eqref{phi1}--\eqref{phi2}, respectively. 

\begin{figure}[t]
\begin{center}
\hspace*{-3mm}
\includegraphics[width=172mm]{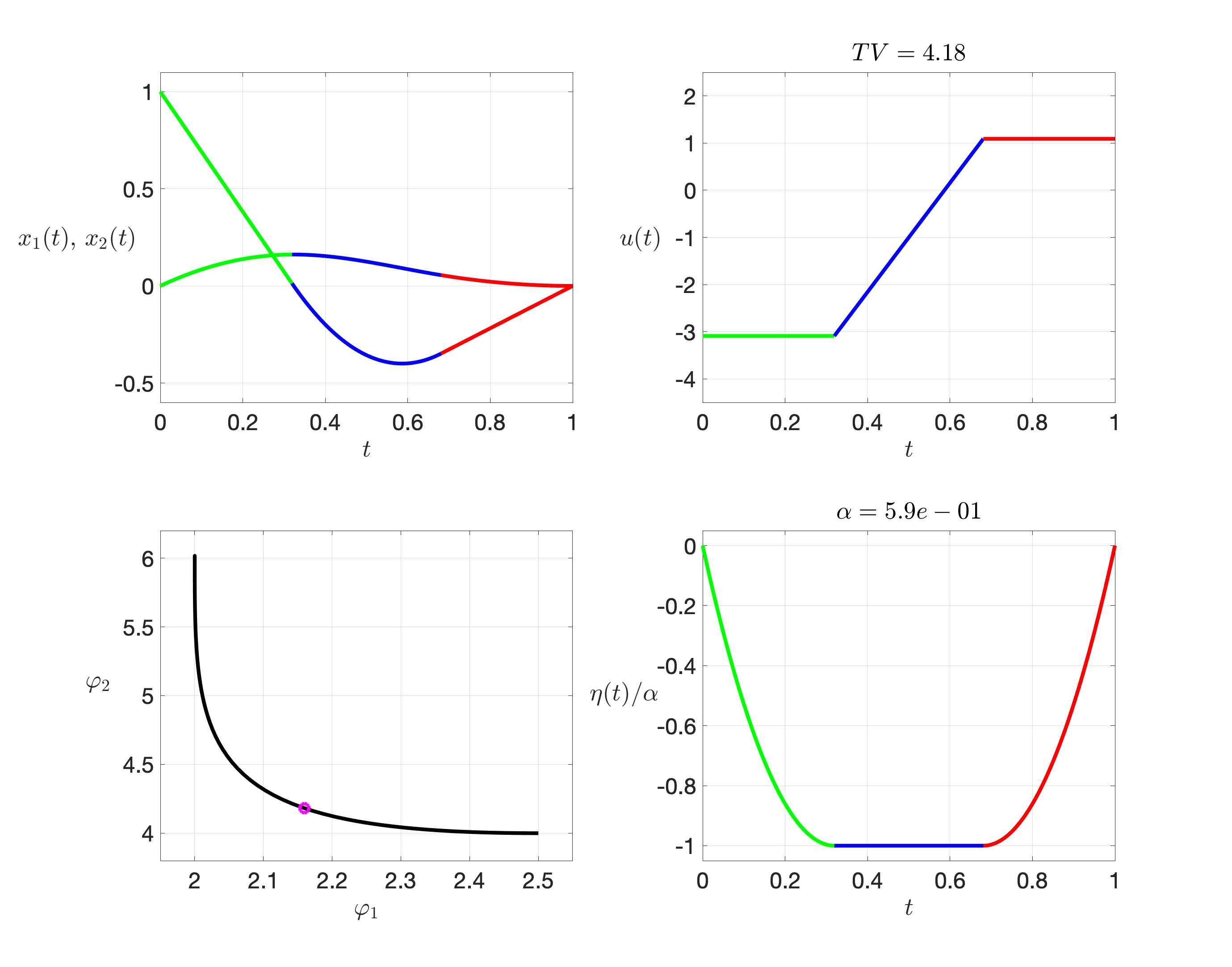}
\end{center}
\caption{\sf A snapshot of the 11th second of the multi-objective solution video {\tt mintotalvar.avi}~\cite{video}, with $\alpha = 0.589$.} 
\label{fig:snapshot}
\end{figure}

Using a rather ``continuous" range of values of $\alpha$, we have generated a movie, by using {\sc Matlab}.  The movie file is called {\tt mintotalvar.avi}, which can be downloaded via the URL in Reference~\cite{video}.  An instance of the movie for $\alpha = 0.589$ is shown in Figure~\ref{fig:snapshot}.  For a large number of values of $\alpha$, the movie depicts/animates the Pareto front (using~\eqref{phi1}--\eqref{phi2}) and the graphs of the control and state variables (using~\eqref{optcont_convex1} and \eqref{x1}--\eqref{x2}), as well as the graph of the adjoint variable $\eta(t)$ divided (or normalized) by $\alpha$ (using~\eqref{eta_convex1}).  The graph of $\eta(t)/\alpha$ in the lower-right corner reconfirms that $u(t)$ is constant when $|\eta(t)| < \alpha$ and $u(t)$ is linear in $t$ when $|\eta(t)| = \alpha$.

As expected, reduction in total control variation is obtained as the value of $\alpha$ is increased, with the trade-off that minimum energy is increased.  Figure~\ref{fig:convex}(b), and the movie, clearly demonstrate that, as $\alpha$ gets larger, the control variable appears to become closer to a piecewise-constant function, switching from the constant level $-3$ to the constant level $1$, resulting in $\tv(u)=4$.  This reconfirms Theorem~\ref{thm:asymp_solution} as well as Remark~\ref{rem:asymptotic}.

Finally, with $\bar{u}_1 = -3$ and $\bar{u}_3 = 1$, the asymptotic expressions for the state variables in \eqref{x1_asymp}--\eqref{x2_asymp} can be rewritten neatly as
\begin{equation}  \label{x1_asymp_neat}
x_1(t) = \left\{\begin{array}{rl}
\ds-\frac{3}{2}\,t^2 + t\,, &\ \ \mbox{if\ \ } 0\le
      t < 1/2\,, \\[3mm] 
\ds\frac{1}{2}\left(t-1\right)^2, &\ \ \mbox{if\ \ } 1/2\le t \le 1\,,
\end{array} \right.
\end{equation}
\begin{equation}  \label{x2_asymp_neat}
x_2(t) = \left\{\begin{array}{rl}
-3\,t + 1\,, &\ \ \mbox{if\ \ } 0\le t < 1/2\,, \\[3mm] 
\ds t - 1\,, &\ \
      \mbox{if\ \ } 1/2\le t \le 1\,.
\end{array} \right.
\end{equation}


\section{Conclusion and Future Work}
\label{sec:conclusion}

We have derived the unique solution to the optimal control problem of simultaneous minimization of energy and total variation in control for the double integrator.  We obtained analytic expressions for the construction of the Pareto front.   We have shown that the unique asymptotic optimal control function, for the minimization of the total variation alone, is piecewise constant with one switching at the midpoint of the time horizon.  We computed the two constant levels of the asymptotic control function analytically.  Subsequently, we have proved that the asymptotic optimal total variation is exactly $2/3$ of the total variation of the minimum-energy control.  These results seem to be the first of their kind in the literature concerning optimal control with minimum total variation, even for a system as simple as the double integrator.

The minimum-energy control problem which we have also considered is a special case of a general linear quadratic control problem.  An approach similar to the one employed in the current paper can be employed for the more general linear quadratic control (or linear quadratic programming) problem where one is additionally concerned with the minimization of total variation, namely the problem
\[
\mbox{(LQPTV) } \left\{\begin{array}{rl}
\min\,& \ \ds
\frac{1}{2}\,\int_{0}^{1} [x(t)^TQ(t)\,x(t)+u(t)^TR(t)u(t)]\,dt + \alpha\tv(u)
\\[5mm]
\mbox {subject to}&\ \dot{x}(t) = A(t)\,x(t) + B(t)\,u(t)\,,\quad\mbox{for all } t\in[0,1]\,, \\[2mm]
 &\ x(0) = x_0\,,\quad  x(1) = x_f\,.
\end{array}\right.
\]
The time horizon in Problem~(LQPTV) has been set to be $[0,1]$, but, without loss of generality, it can be taken to be any interval $[t_0,t_f]$, with $t_0$ and $t_f$ specified. The state variable vector $x(t)\in\dR^n$ and the control variable vector $u(t)\in\dR^m$. The time-varying matrices $A:[0,1]\to \dR^{n\times n}$ and $B:[0,1]\to\dR^{n\times m}$ are continuous, $Q:[0,1]\to \dR^{n\times n}$ is symmetric positive definite and continuous in $t$, and $R:[0,1]\to \dR^{m\times m}$ is positive definite and continuous in $t$. The initial and terminal states are specified as $x_0$ and $x_f$, respectively.  Since there are more than just one control variable, i.e., $u(t) = (\bar{u}_1(t),\ldots,u_m(t))\in\dR^m$, the total variation in \eqref{totalvar2} can be generalized for this case as
\begin{equation} \label{totalvar3}
\tv(u) := \int_0^1 \left(|\dot{u}_1(t)| + \ldots + |\dot{u}_m(t)|\right)\,dt\,. 
\end{equation}
It should be noted that the problem we have studied in the current paper fits into the above problem description (LQPTV) with $n=2$, $m=1$, $Q=0$ and $R=1$, and the appropriate constant system and control matrices $A$ and $B$.

The general linear quadratic problem is a convex problem, so the weighted-sum scalarization can still be used (see \cite{BonKay2010,KayMau2014}) when it is combined with the minimization of total variation.  However, for a generalization to nonconvex problems, a scalarization different from the weighted-sum scalarization needs to be considered.  This requires specialized numerical techniques in obtaining a solution---see \cite{KayMau2014} and the pertaining discussion therein for problems which also have constraints on the state and control variables.

\section*{Acknowledgments}

The author would like to offer his warm thanks to two anonymous reviewers, whose comments and suggestions improved the paper.

\end{document}